\documentclass[12pt,a4paper]{amsart}
\usepackage{amsmath,amssymb,amsthm}
\usepackage{enumitem}
\usepackage{ifthen,verbatim}
\usepackage{mathrsfs}
\usepackage{color}
\usepackage{url}
\usepackage[english]{babel}
\usepackage{here}
\usepackage[T1]{fontenc}
\usepackage{floatflt,graphicx,graphics}
\usepackage{a4wide}
\usepackage{cite}

\numberwithin{equation}{section}
\nonstopmode
\theoremstyle{plain}
\newtheorem{thm}[equation]{Theorem}
\newtheorem{cor}[equation]{Corollary}

\theoremstyle{definition}

\theoremstyle{remark}

\newtheorem{nonsec}[equation]{}

\numberwithin{equation}{section}

{\qed\bigskip}

\DeclareMathOperator{\am}{am}
\font\fFt=eusm10 
\font\fFa=eusm7 
\font\fFp=eusm5 
\def\K{\mathchoice
{\hbox{\,\fFt K}}
{\hbox{\,\fFt K}}
{\hbox{\,\fFa K}}
{\hbox{\,\fFp K}}}

\newcommand{\cn}{\mathrm{cn}}
\newcommand{\dn}{\mathrm{dn}}
\newcommand{\sn}{\mathrm{sn}}
\newcommand{\D}{{\mathbb D}}
\renewcommand{\i}{\mathrm{i}}
\renewcommand{\Re}{{\,\operatorname{Re}\,}}
\renewcommand{\Im}{{\,\operatorname{Im}\,}}

\newcounter{minutes}
\divide\time by 60
\newcounter{hours}
\multiply\time by 60
\addtocounter{minutes}{-\time}

\begin{document}

\title[]{LEVEL SETS OF POTENTIAL FUNCTIONS\\
BISECTING UNBOUNDED QUADRILATERALS}
\author[M.M.S.~Nasser]{Mohamed M.~S.~Nasser}
\address{Program of Mathematics, Department of Mathematics, Statistics and Physics, College of Arts and Sciences, Qatar University, Doha, Qatar}
\email{mms.nasser@qu.edu.qa}
\author[S.~Nasyrov]{Semen Nasyrov}
\address{Institute of Mathematics and Mechanics, Kazan Federal University, 420008 Kazan, Russia}
\email{semen.nasyrov@yandex.ru}
\author[M. Vuorinen]{Matti Vuorinen}
\address{Department of Mathematics and Statistics, University of Turku, FI-20014 Turku, Finland}
\email{vuorinen@utu.fi}

\date{}


\begin{abstract} We study the mixed Dirichlet-Neumann problem for the Laplace equation
in the complement of a bounded convex polygonal quadrilateral in the extended complex plane.
The Dirichlet\,/\,Neumann conditions at opposite pairs of sides are $\{0,1\}$ and
$\{0,0\},$  resp. The solution to this problem is a harmonic function in the unbounded complement of the polygon known as the \emph{potential function} of the quadrilateral. We compute the values of the potential function including its value at infinity.
\end{abstract}



\maketitle

\def\thefootnote{}
\footnotetext{\texttt{{\tiny File:~\jobname .tex, printed: \number\year-%
\number\month-\number\day, \thehours.\ifnum\theminutes<10{0}\fi\theminutes}}}
\makeatletter

\makeatother

\section{Introduction}

A \emph{quadrilateral} in the extended complex plane is a Jordan domain $D$ together with
an ordered quadruple $\{z_1,z_2,z_3,z_4\}$ of points on its boundary. The points are
ordered in such a way that when traversing the boundary, the points occur in the order
of the indices and the domain $D$ is on the left-side. Such a quadrilateral is denoted
$Q=(D; z_1,z_2,z_3,z_4)$. By the Riemann mapping theorem, the simply connected domain
$D$ can be mapped by a conformal mapping
onto a rectangle and by Caratheodory's theorem a conformal mapping between two
Jordan domains extends to a homeomorphism between their closures.
Therefore there exists a conformal mapping $f$ of $D$ which defines
a homeomorphism
\begin{equation}\label{eq:f}
f:\overline{D} \to [0,1] \times [0, h],\quad h>0,
\end{equation}
such that  \cite[ p. 52]{papa}
\begin{equation}\label{eq:fcd}
f(z_1)= 0,\ f(z_2)=1,\ f(z_3)=1+\i h,\ f(z_4)= \i h.
\end{equation}
The parameter $h= {\rm Mod}(Q)$ is an important domain characteristic, it is the
\emph{modulus of the quadrilateral,} widely studied in geometric function theory and related areas of geometric analysis~\cite{ahlfors,d,gm,hkv}. It is a basic property that the modulus
is a conformal invariant and numerous applications of the modulus depend on this fact.

Here, we will analyse
polygonal quadrilaterals. Note that  after relabeling the points $z_j$ in the opposite order the domain  complementary  to $D$ together with the relabelled points forms another quadrilateral. If $D$ is bounded, we call the modulus of $D$ the \emph{interior modulus} and the modulus of the complementary
domain the \emph{exterior modulus}. In the case of a bounded convex polygonal quadrilateral one can apply the Schwarz-Christoffel transformation and give an explicit formula for the interior modulus in terms of the vertices of the quadrilateral and special functions as shown in \cite[Section 2]{hvv}. The case of the exterior modulus is similar, but much more involved, as shown very recently in \cite{nsv}.

There are many methods designed for numerical conformal mapping,
see the recent short surveys \cite[pp.14-16]{papa}, \cite[pp.8-12]{ky},
and the long survey  \cite{w}.
The method best known for Schwarz-Christoffel problems  due to N.~Trefethen
is implemented in the MATLAB Toolbox of T. Driscoll \cite{dt}.
Two numerical methods have proven successful not only for
computation of the moduli of polygonal quadrilaterals but for a wide class of
domains with curvilinear boundaries which occur in
applications. These two methods are the $hp$-FEM method \cite{hrv1,hrv2}
and the boundary integral equation method~\cite{nasser,nrv}.
The $hp$-FEM method makes use of an alternative equivalent
definition of the modulus of a quadrilateral in terms of a mixed
Dirichlet-Neumann boundary value problem \cite[Thm 4.5, p. 63]{ahlfors},
\cite{d}. Consider a quadrilateral $Q = (D; z_1, z_2, z_3, z_4)$ where
$D$ is a Jordan domain on the Riemann sphere containing $\infty$.
Let $u$ be a harmonic function in $D$ solving the following mixed
Dirichlet-Neumann boundary value problem:
\begin{equation} \label{eq:dirichlet}
\left\{\begin{matrix}
    \Delta u& =&\ 0,& \ \text{on}\ &{\ D,} \\
    u& =&\ 0,&  \ \text{on}\ &{\partial Q_1 = (z_4,z_1),}\\
    u& =&\ 1, & \ \text{on}\ &{\partial Q_3 = (z_2,z_3),}\\
    \partial u/\partial n&  =&\ 0,& \ \text{on}\ &{\partial Q_2 = (z_1,z_2),}\\
    \partial u/\partial n&  =&\ 0,& \ \text{on}\ &{\partial Q_4 = (z_3,z_4).}\\
\end{matrix}\right.
\end{equation}

It is well-known that this harmonic function $u: D \to (0,1)$ is unique.
It is called the \emph{potential function} of the quadrilateral $Q$.
Because harmonic functions satisfy the maximum principle, the level
sets $\{z \in D: u(z) = t\}$, $0<t<1$, cannot be compact subsets of $D.$
Moreover, the potential function $u$ has a limit at infinity:
$$
u(\infty) = \lim_{|z|\to \infty} u(z).
$$
In particular, only one of the level curves of the potential function passes
through $\infty$ and this level set bisects the domain $D$ into
two unbounded parts.

We note that the potential function $u$ is the real part of
the function $f$ given by \eqref{eq:f}. Therefore, we have
$f=u+iv$ where $v$ is the harmonic function conjugate to $u$.
Then the function $g=-i(f-1)/h$ maps conformally $D$ onto the
rectangle $[0,1]\times[0,1/h]$ and $v/h=\Re g$. Thus, $v/h$ is
the potential function for the conjugate quadrilateral
$Q^*=(D; z_2,z_3,z_4,z_1)$. Finding the values of
$f(\infty)=u(\infty)+iv(\infty)$ and $\lim_{z\to\infty}[z(f(z)-f(\infty))]$
also allows us to calculate the Robin constant of $D$ which is
closely connected with the logarithmic capacity of its complement
$K=D^c$, the transfinite diameter and the Chebyshev constant of $K$
(see, e.g., \cite[ch.7]{gol}, \cite{LSN}).
We also note the connections of these concepts with the reduced modulus
of $D$ at infinity \cite{d,gm}.

The main result of this note is Theorem \ref{inftyMain} based on
the recent work \cite{nsv}; this theorem yields a formula
for $u(\infty)$.
The novelty of this result lies in the fact that the formula is explicit,
expressed  in terms of the angles of the given quadrilateral and the
well-known special functions.

Previously, the problem of computing the exterior modulus of a bounded
convex polygonal quadrilateral was studied in \cite{hrv2}; 
the idea was to use the inversion transformation $z  \mapsto z/|z|^2$ and
thereby to reduce the problem to the case of a bounded domain bounded by
circular arcs and to the numerical solution of a mixed Dirichlet-Neumann
problem in such a domain. This method does give a numerical approximation
of the value of the potential function at infinity, but the formula
of Theorem \ref{inftyMain} is analytic.

We also use two independent numerical methods to illustrate our results,
to provide a short table of the values of $u(\infty)$ and graphics for the
level curves of the potential function $u$ for two specific bounded
convex polygonal quadrilaterals. The first method is a Mathematica program
adopted from \cite{nsv} for the present purpose.
The second method is based on using the MATLAB toolbox PlgCirMap from~\cite{pcm} to compute the conformal mapping $f$ in~\eqref{eq:f}--\eqref{eq:fcd}. With the help of this toolbox, the modulus $h= {\rm Mod}(Q)$ will be computed in a similar way as in~\cite{nasser}.
For the purpose of visual comparison, we give graphics of level sets of the potential function produced by both methods. The tabular data obtained show that the two methods agree with high precision, at least with 9 decimal places.

\section{Notation}

Conformal mappings are often given by special functions \cite{af,a,ky}.
Therefore, various special functions
are recurrent in the study of moduli of quadrilaterals. We need here complete elliptic integrals,
$\K(r)$ and $\K'(r)$  of the first kind defined as \cite{avv}
\begin{equation} \label{eq:ellipK}
\K(r)=\int^1_0 \frac{dx}{\sqrt{(1-x^2)(1-r^2x^2)}}\,,
\quad \K'(r)=\K(r'), \quad r'=\sqrt{1-r^2}\,.
\end{equation}
Equivalently, one can use the Gaussian hypergeometric function ${}_2 F_1$ to define
$\K:$
$$
\K(r) = {\frac{\pi}{2}}\,{}_2 F_1\bigl( {\textstyle \frac{1}{2},\frac{1}{2}};1; r^2\bigr),\quad 0 <r<1.
$$
The incomplete elliptic integral is defined by
\begin{equation} \label{eq:ellipF}
F(r,\varphi)= \int^{\varphi}_0 \frac{d\theta}{\sqrt{1-r^2 \sin^2 \theta}}=
\int^z_0 \frac{dt}{\sqrt{(1-t^2)(1-r^2t^2)}} \quad (0\le r<1,\ z = \sin \varphi).
\end{equation}
Clearly, $\K(r) = F(r, \pi/2).$
For more information on the complete
elliptic integral $\K$
and other special functions, see
\cite{avv,bateman}.

We will also need the Jacobi sine function $\sn(z,\lambda)$.
By definition (see, e.g. \cite[Ch.V, \S 24]{a}),
$\sn(z,\lambda)=\sin \am(z,\lambda)$  where $\am(z,\lambda)$
is the Jacobi amplitude, i.e. the function $\varphi=\am(z,\lambda)$
inverse to the incomplete elliptic integral of the
first kind
$$
z=\int_0^\varphi\frac{d \theta}{\sqrt{1-\lambda^2\sin^2\theta}}\,.
$$
Elliptic integrals and elliptic functions occur for instance when
the upper half plane is conformally mapped onto a rectangle \cite[p. 358]{af}.

\section{Numerical computation using PlgCirMap MATLAB toolbox}

The modulus of quadrilaterals and the conformal mappings for polygonal domains can be computed by the numerical method presented in~\cite{nasser}. The method is based on using the boundary integral equation with the generalized Neumann kernel.
In this paper, instead, we will use the PlgCirMap MATLAB toolbox~\cite{pcm} which is also based on using the boundary integral equation with generalized Neumann type kernel to compute the conformal mapping from polygonal domains onto circular domains and its inverse. It provides a flexible computational method for several problems in the area of applied and computational complex analysis.
The toolbox will be used to compute the conformal mapping $f$ in~\eqref{eq:f}, the modulus $h= {\rm Mod}(D; z_1,z_2,z_3,z_4)$, as well as the potential function
$u(z)$ for $z\in D$.
In the numerical computation below, we choose the number of discretization points on each side of the polygons to be $n=2^{13}$. Thus, the total number of discretization points on the boundary of the polygonal quadrilateral is $2^{15}$.

\begin{figure}[H] \centering
\includegraphics[width=8cm]{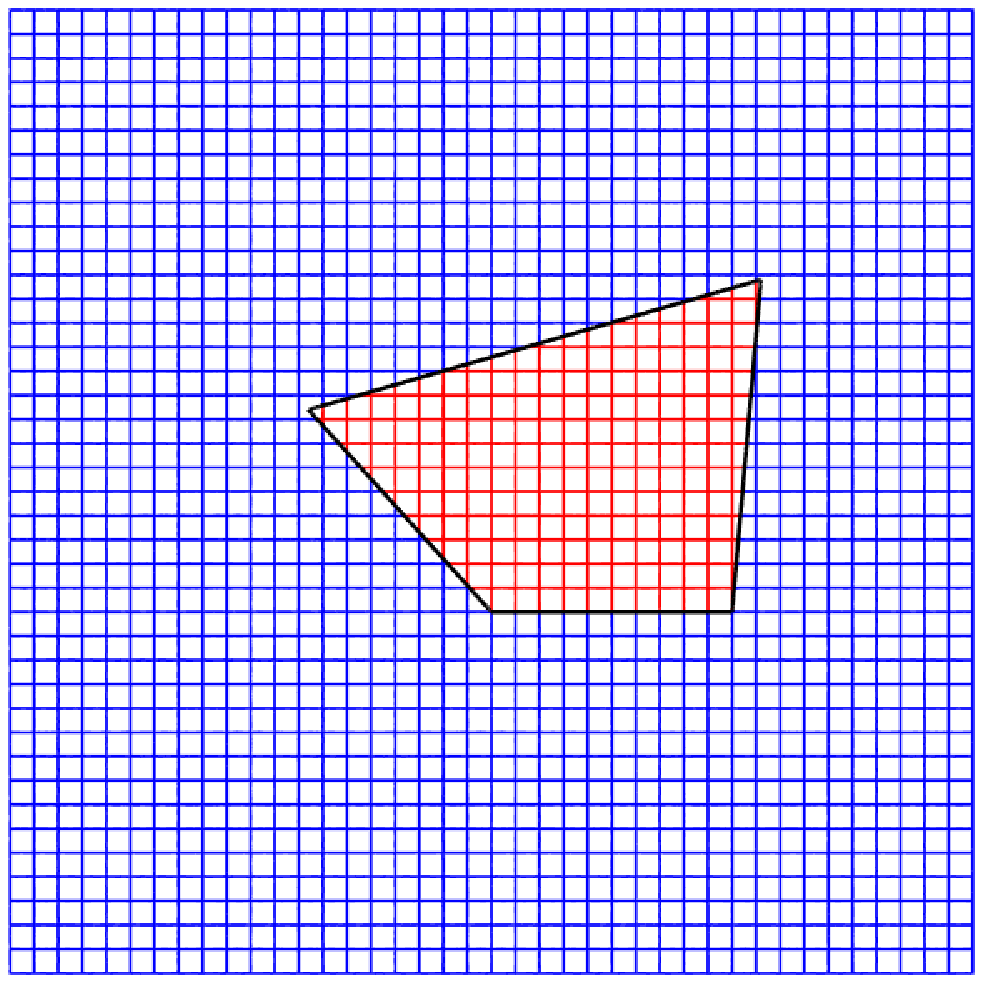}
\hfill
\includegraphics[width=8cm]{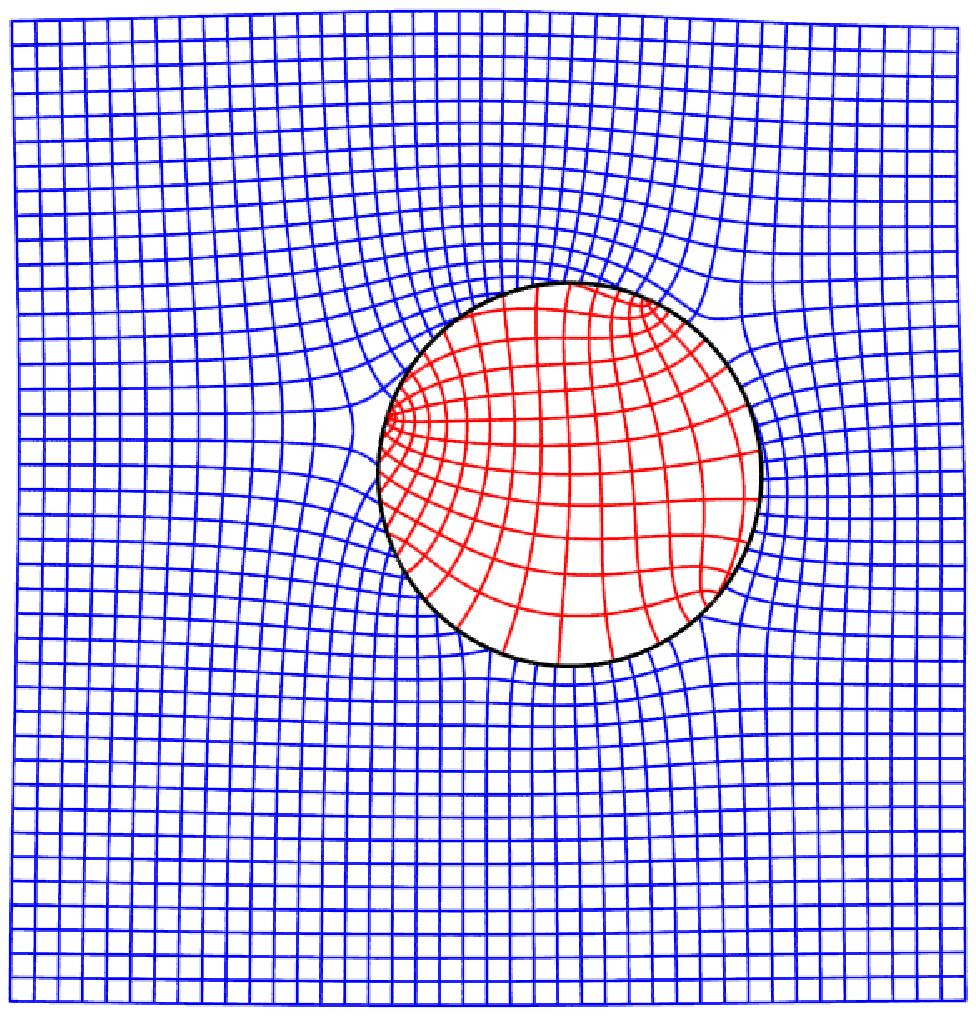}
\caption{A coordinate grid (left) and its image under two conformal
maps $\Phi_1$ and $\Phi_2$, defined in the exterior and in the interior of the polygonal quadrilateral of Example \ref{ep1}, resp.
The exterior map $\Phi_1$ is normalized by assuming that $\Phi_1(z)=z+O(1/z)$ near infinity.
The interior map $\Phi_2$ is normalized by assuming that $\Phi_2(\hat\alpha)=0$ and $\Phi'_2(\hat\alpha)>0$ where $\hat\alpha$ is a given interior point to the quadrilateral.
}
\label{fig0}
\end{figure}

\begin{nonsec}{\bf The exterior and interior of a polygonal domain.}
We demonstrate here the toolbox, by constructing for a convex
bounded polygonal quadrilateral two conformal mappings: the first mapping maps
the interior of a polygonal line $L$ onto the unit disk and the other
one mapping the exterior of $L$
onto the exterior of the unit disk. The same method would also work  for much
more general polygonal domains, such as Koch's snowflake domain
considered in~\cite[Figure~2.1]{wk}. In Figure \ref{fig0} we
show how the rectangular coordinate grid, on the left side of the figure
is transformed onto a curvilinear grid. Note  on the right side of the figure,
the discontinuous behavior of the two image grids
at the points of the unit circle.
The images of the coordinate grid shows how the
exterior and the interior of the quadrilateral
are conformally different. For both the conformal mappings,
the strongest distortion
occurs close to the vertices.  By Caratheodory's
theorem, each of the conformal mappings has a homeomorphic
extension to the closure of the quadrilateral, but
the values of the exterior and interior
conformal maps do not agree on the boundary of the quadrilateral.
This polygonal domain will be studied below in Example~\ref{ep1}.
\end{nonsec}

\begin{nonsec}{\bf Computing the quadrilateral modulus.}
To compute the modulus $h= {\rm Mod}(D; z_1,z_2,z_3,z_4)$, we first compute the unique conformal mapping
$\Phi_1(z)$ from the unbounded domain $D$ onto the exterior of the
unit disk with the normalization that $\Phi_1(z)=z+O(1/z)$ near infinity.
This mapping function is computed by calling
\[
{\tt Phi1   =  plgcirmap(ver,inf)}
\]
where ${\tt ver\{1\}=[z_1,z_2,z_3,z_4]}$. Then, the conformal mapping
\[
\hat z=\Phi(z)=\frac{1}{\Phi_1(z)}
\]
maps the unbounded domain $D$ onto the unit disk
$\D=\{\hat z\,:\,|\hat z|<1\}$. The boundary values
$\hat\eta(t)=\Phi(\eta(t))$
 form a parametrization of the unit circle $\partial\D$
 and can be computed by calling
\[
\mbox{\tt het = Phi1.zet}.
\]
Hence, the vertices $z_k$, $k=1,2,3,4$, of the polygon
domain $D$ are mapped onto the points $\hat z_k$, $k=1,2,3,4$,
on the unit circle where
\[
\hat z_k={\tt het((k-1)n/4+1)}, \quad k=1,2,3,4.
\]
Thus, by the conformal invariance of the modulus, the modulus of the quadrilateral $(D; z_1,z_2,z_3,z_4)$ is given by~\cite{hkv,nasser,papa}
\begin{equation}\label{eq:mod}
h={\rm Mod}(D;z_1,z_2,z_3,z_4) =
{\rm Mod}(\D;\hat z_1,\hat z_2,\hat z_3,\hat z_4) =
 \frac{2}{\pi}\,\mu\left(1/\sqrt{k}\right)
\end{equation}
where
\begin{equation}\label{eq:ration}
\mu(r) =  \frac{\pi \, \K(r')}{2 \K(r)}\,, \quad
k=|\hat z_1,\hat z_2,\hat z_3,\hat z_4| = \frac{|\hat z_1-\hat z_3|\,|\hat z_2-\hat z_4|}{|\hat z_1-\hat z_2|\,|\hat z_2-\hat z_4|}\,.
\end{equation}

\end{nonsec}

\begin{figure}[ht]
\centerline{
\includegraphics[width=11cm,trim=0 1.0cm 0 0,clip]{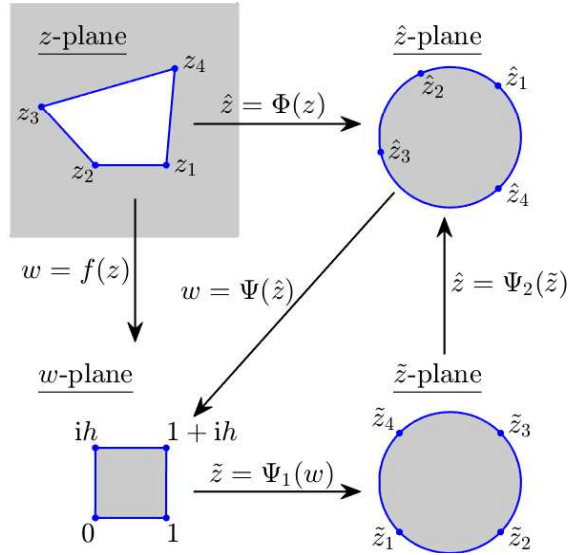}}
\caption{The numerical computation of the conformal mapping $w=f(z)$ using the PlgCirMap MATLAB toolbox.}\label{fig:map}
\end{figure}

\begin{nonsec}{\bf Computing the conformal mapping $f$ in~\eqref{eq:f}.}
Here, we describe how the toolbox can be used to compute the conformal mapping $f$ from $D$ onto the
rectangle $R=\{w\,:\, 0<\Re w<1,\,0<\Im w<h\}$ in~\eqref{eq:f} that satisfies the condition~\eqref{eq:fcd}. The method is summarized in Figure~\ref{fig:map}.

By the definition of $h={\rm Mod}(\D;\hat z_1,\hat z_2,\hat z_3,\hat z_4)$
\cite[p.52]{papa}, there exists a conformal mapping
\[
w=\Psi(\hat z)
\]
from the unit disk $\D$ onto the rectangle $R$ such that
\[
\Psi(\hat z_1)=0,\;\Psi(\hat z_2)=1,\;\Psi(\hat z_3)=1+\i h,\;
\Psi(\hat z_4)=\i h.
\]
To compute such a conformal mapping $\Psi$, we first compute
the unique conformal mapping
\[
\tilde z=\Psi_1(w)
\]
from the domain $R$ onto the unit disk $\D$ with the normalization
\begin{equation}\label{eq:Phi-cond-b}
\Psi_1(\hat\alpha)=0, \quad \Psi_1'(\hat\alpha)>0
\end{equation}
where $\hat\alpha$ is an auxiliary point in $R$, say $\hat\alpha=(1+\i h)/2$.
 This conformal mapping $\Psi_1$ can be computed by the MATLAB toolbox
 PlgCirMap by calling
\[
{\tt Psi1   =  plgcirmap(vers,halpha)}
\]
where ${\tt vers\{1\}}=[0,1,1+\i h,\i h]$ is the vector of
vertices of $\partial R$ and ${\tt halpha}=\hat\alpha$.
 The mapping function $\tilde z=\Psi_1(w)$ maps the
 vertices $0,1,1+\i h,\i h$ of $\partial R$ onto four points
 $\tilde z_1,\tilde z_2,\tilde z_3,\tilde z_4\in\partial\D$ where
\[
\tilde z_k={\tt tzet((k-1)n/4+1)}, \quad k=1,2,3,4,
\]
and {\tt tzet   =  Psi1.zet}.
These points are in general different from the points
$\hat z_1,\hat z_2,\hat z_3,\hat z_4$. Let
\[
\hat z=\Psi_2(\tilde z)= \hat z_3+
\frac{(\hat z_3-\hat z_1)(\hat z_2-\hat z_3)(\tilde z_2-\tilde z_1)(\tilde z-\tilde z_3)}{(\hat z_2-\hat z_1)(\tilde z_2-\tilde z_3)(\tilde z-\tilde z_1)-(\hat z_2-\hat z_3)(\tilde z_2-\tilde z_1)(\tilde z-\tilde z_3)}\,,
\]
then $\Psi_2$ maps the unit disk $\D$ onto itself such
that $\Psi_2(\tilde z_1)=\hat z_1$,
$\Psi_2(\tilde z_2)=\hat z_2$, and $\Psi_2(\tilde z_3)=\hat z_3$.
Thus, the function
\[
w=(\Psi_1^{-1}\circ\Psi_2^{-1})(\hat z)
\]
maps the unit disk $\D$ onto the rectangle $R$ and takes
the three points $\hat z_1$, $\hat z_2$, $\hat z_3$ to the three points $0$, $1$, $1+\i h$, respectively.
Since the function $\Psi$ is also a conformal mapping from the unit disk $\D$ onto the rectangle $R$ and maps the three points $\hat z_1,\hat z_2,\hat z_3$ to the three points $0,1,1+\i h$, respectively, then we have
\[
\Psi = \Psi_1^{-1}\circ\Psi_2^{-1}.
\]
This is due to the fact of the uniqueness of the conformal mapping that maps the unit disk $\D$ onto the domain $R$ and maps three points on $\partial D$ to three points on $\partial R$ when $h$ is fixed.

Thus, the function
\[
w=f(z)=(\Psi\circ\Phi)(z)=(\Psi_1^{-1}\circ\Psi_2^{-1}\circ\Phi)(z)
\]
is the required unique conformal mapping satisfying~\eqref{eq:f} and~\eqref{eq:fcd}.

\end{nonsec}

\begin{nonsec}{\bf Computing the values of the potential function $u$.}
Finally, we describe how to use the toolbox to compute the values of the mapping function $w=f(z)=\Psi_1^{-1}(\Psi_2^{-1}(\Phi(z)))$ and then the values of the potential function $u(z)$ for interior points $z\in D$.

First, for $z\in D$, the values of the mapping function $\hat z=\Phi(z)$ can be computed using the MATLAB function \verb|evalu| from the PlgCirMap MATLAB toolbox by calling
\[
{\tt \hat z = 1./evalu(Phi1,z,{'d'})}
\]
For the M\"obius transform, $\Psi_2^{-1}:\D\to \D$, it can be computed by the explicit formula
\[
\tilde z=\Psi_2^{-1}(\hat z)= \tilde z_3+
\frac{(\tilde z_3-\tilde z_1)(\tilde z_2-\tilde z_3)(\hat z_2-\hat z_1)(\hat z-\hat z_3)}{(\tilde z_2-\tilde z_1)(\hat z_2-\hat z_3)(\hat z-\hat z_1)-(\tilde z_2-\tilde z_3)(\hat z_2-\hat z_1)(\hat z-\hat z_3)}\,.
\]
Then, the values of the inverse mapping function $w=\Psi_1^{-1}(\tilde z)$ can be computed by the MATLAB function \verb|evalu| by calling
\[
{\tt w=evalu(Psi1,\tilde z,{'v'})}
\]

By computing the values of the function $w=f(z)$, we find the values of the potential function $u$ through
\[
u(z)=\Re f(z), \quad z\in D.
\]
Since $\Phi(\infty)=0$, we have
\[
u(\infty) = \Re f(\infty)= \Psi_1^{-1}(\Psi_2^{-1}(\Phi(\infty))) = \Psi_1^{-1}(\Psi_2^{-1}(0)).
\]

\end{nonsec}

\section{The main result and its proof}

We will apply \cite{nsv} to give  a formula for the value of the potential function
at infinity, $u(\infty),$ for a quadrilateral $Q=(D; z_1,z_2,z_3,z_4)$ in the case where $D$ is an unbounded domain with a boundary which is a convex polygonal line with vertices $z_1$, $z_2$, $z_3$, $z_4$,
and interior angles of $D$ equal to
\begin{equation}\label{angles}
(1-\alpha)\pi,\;(1-\beta)\pi,\; (1-\gamma)\pi,\; (1-\delta)\pi, \quad  \alpha,
\beta, \gamma, \delta\in(0,1), \quad \alpha+\beta+\gamma+\delta=2.
\end{equation}
This formula for the potential function, given in Corollary \ref{infinit}, is our
main result; it is expressed in terms of special functions.

Let $t=1/r^2,$ where $r\in(0,1),$ be the unique root of the equation
\begin{equation}\label{rk}
\text{Mod}(Q)={\K}(r')/{\K}(r), \quad
r'=\sqrt{1-r^2};
\end{equation}
here $\text{Mod}(Q)$ is the conformal modulus of $Q$.
Denote $\mathcal{E}=\alpha+\beta+\gamma-1$,
$$
\mathcal{A}=2 (\mathcal{E}-1)^2,\quad
\mathcal{B}=(\mathcal{E}-1)[4-3(\alpha+\gamma)+(4-3(\alpha+\beta))t],
$$
$$
\mathcal{C}=2-3(\alpha+\gamma)+(\alpha+\gamma)^2+2(3-5\alpha-2 \beta-2
\gamma+2\alpha^2+2\alpha \beta+2 \alpha \gamma+\beta \gamma)t$$
$$+(2-3(\alpha+\beta)+(\alpha+\beta)^2)t^2,
\quad  \mathcal{D}=(1-\alpha)(\alpha+\gamma-1+(\alpha+\beta-1)t)t,
$$
\begin{equation*}\label{rho}
\rho(x)=\frac{\alpha t x } {(1-\mathcal{E})x+\mathcal{E}(t+1)-\gamma t-\beta}\,.
\end{equation*}

In \cite{nsv} the following result is actually proved.

\begin{thm} \label{inftyMain}
The conformal mapping $f$ of the upper half plane
onto $D$ is given
by the generalized Schwarz-Christoffel formula \cite[ Section  5.6, formula  (5.6.3b)]{af}
\begin{equation}\label{F}
f(z)=C\int_0^z
\frac{\zeta^\alpha(\zeta-1)^\beta(\zeta-t)^\gamma}{(\zeta-z_0)^2(\zeta-\overline{z}_0)^2}\,d\zeta+z_1,
\end{equation}
$$
C=(z_2-z_1)\biggl(\int_0^1
\frac{\zeta^\alpha(\zeta-1)^\beta(\zeta-t)^\gamma}{(\zeta-z_0)^2(\zeta-\overline{z}_0)^2}\,d\zeta\biggr)^{-1}.
$$
Here $z_0=x_0+iy_0$ is a pole of $f$ satisfying the following equation:
\begin{equation}\label{eqz0}
\frac{\alpha}{z_0}+\frac{\beta}{z_0-1}+\frac{\gamma}{z_0-t}=\frac{1}{\i y_0}\,.
\end{equation}
The equation \eqref{eqz0} has a unique solution $z_0$ in the upper half plane.
Moreover,
$x_0$ is the unique solution of the
cubic equation
\begin{equation*}\label{4a}
\mathcal{A}x^3+\mathcal{B}x^2+\mathcal{C}x+\mathcal{D}=0,
\end{equation*}
satisfying the inequality $x^2<\rho(x)$, and
$y_0=\sqrt{\rho(x_0)-x_0^2}$.
\end{thm}

Now we will give a method to find $u(\infty)$ with the help of elliptic integrals of the first kind.

Let, as above, $D$ be an unbounded domain with boundary which is a convex polygonal line with vertices
$z_1$, $z_2$, $z_3$, and $z_4$, and interior angles of $D$ satisfying \eqref{angles}.  Let $r\in(0,1)$ be a unique solution of \eqref{rk} and $t=1/r^2$.

\begin{cor}\label{infinit}
For the unbounded polygonal quadrilateral $(D;z_1,z_2,z_3,z_4)$, the potential function \eqref{eq:dirichlet}
is equal to $\Re \psi
\circ f^{-1}$ where $f$ is given by \eqref{F},
$$
\psi(z)=\frac{F(r,\arcsin\sqrt{z})}{\K(r)}\,,
$$
$F(r,\varphi)$ and $\K(r)$ are elliptic integrals given by \eqref{eq:ellipF} and \eqref{eq:ellipK}.
Moreover,
\begin{equation}\label{uinf}
u(\infty)=\Re\, \psi(z_0)
\end{equation}
where $z_0$ is the unique root of \eqref{eqz0} lying in the upper half plane.
\end{cor}

\begin{proof} Actually, $\psi$ maps conformally the upper half plane onto the rectangle $[0,1]\times[0,h]$ and this immediately implies the statement of Corollary~\ref{infinit}.
\end{proof}

\begin{figure}[h]
\begin{minipage}[h]{0.49\linewidth}
\center{\includegraphics[width=0.8\linewidth]{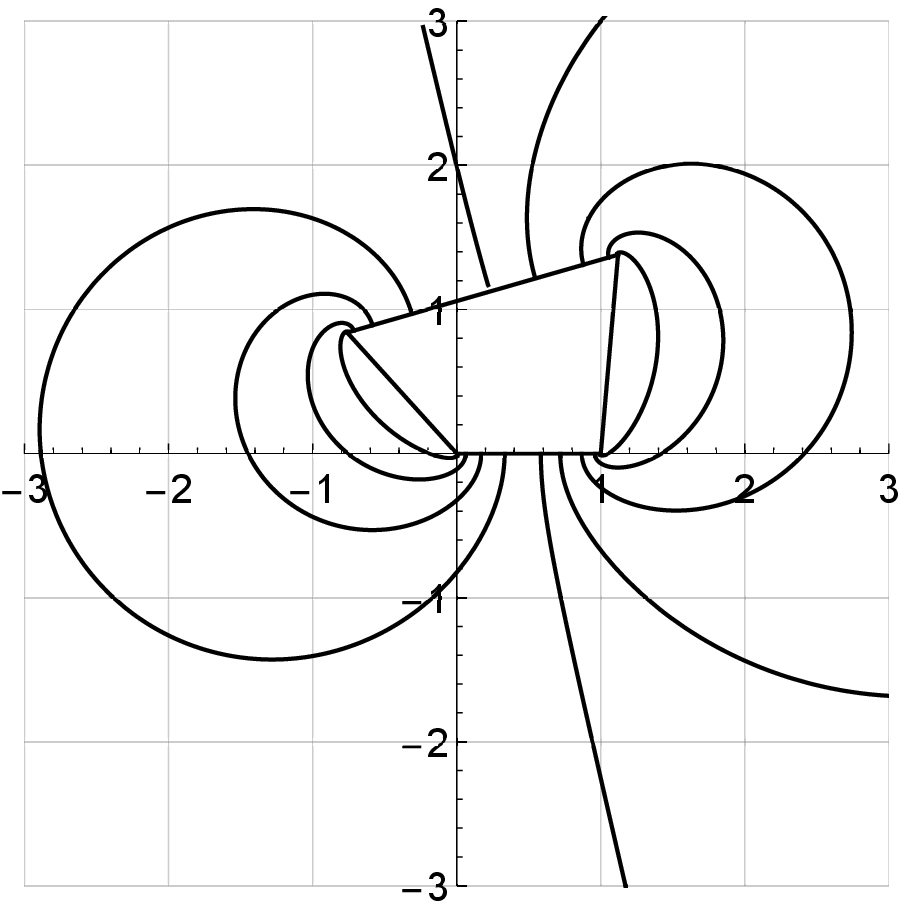}}
\end{minipage}
\hfill
\begin{minipage}[h]{0.49\linewidth}\vspace{3mm}
\center{\includegraphics[width=0.83\linewidth]{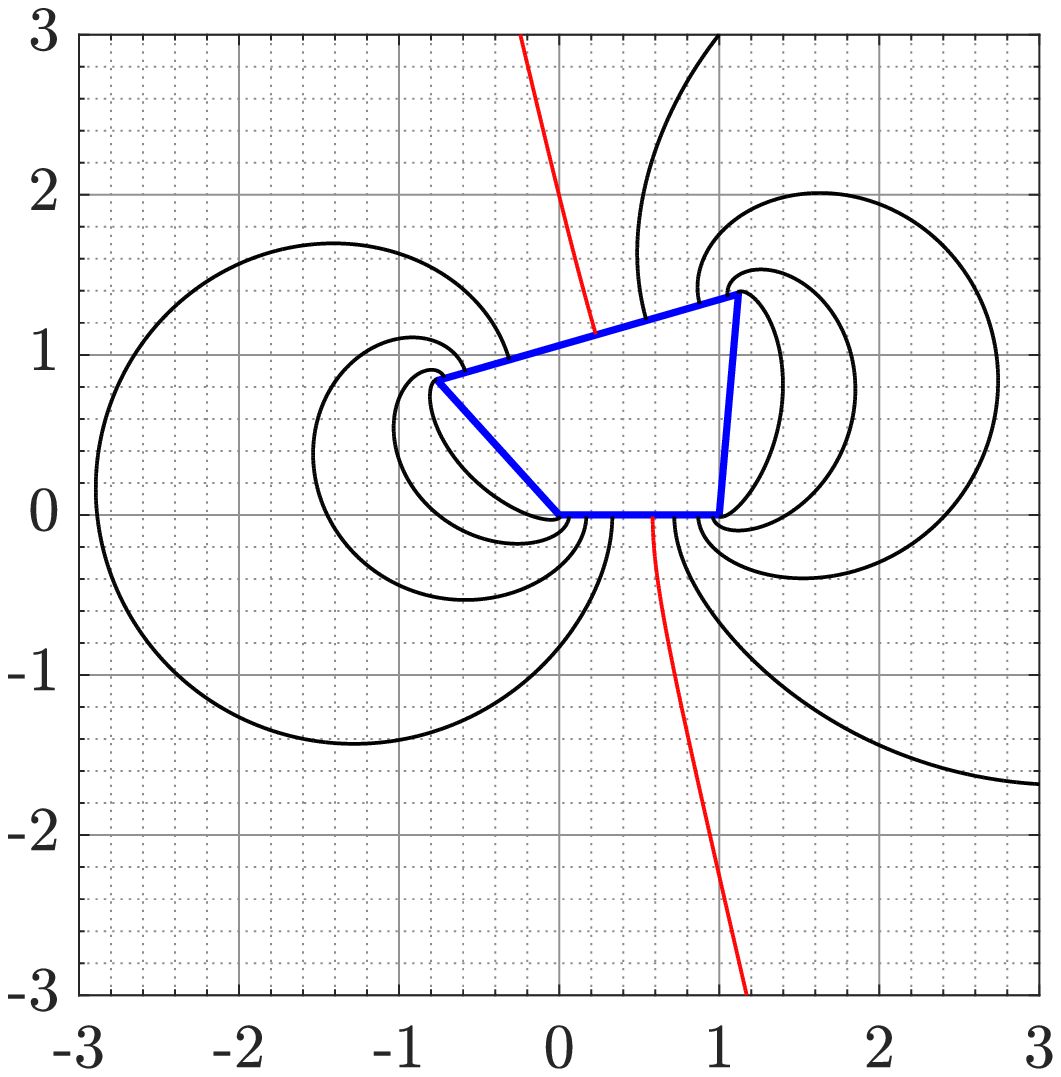}}
\end{minipage}
\caption{The exterior polygon \ref{ep1} and the level curves of its potential function. The figure on the left/right side is produced with the Mathematica/MATLAB codes.}
\label{extp}
\end{figure}

\begin{nonsec}{\bf Some level sets of the exterior modulus.}\label{ep1}
 The polygon shown in Fig.\ref{extp} has vertices
$1, 0, -19/25 + \i\, 21/25, 28/25 + \i\, 69/50.$
 Some level sets of the exterior modulus are shown, corresponding to values:
$0.1$, $0.2$, $0.3$, $0.4$, $u(\infty)=0.471813...$, $0.6$, $0.7$, $0.8$ and $0.9$.
 Compare with \cite[Fig. 3.5 (b)]{hrv2}.
\end{nonsec}



\begin{figure}[h]
\begin{minipage}[h]{0.49\linewidth}
\center{\includegraphics[width=0.7\linewidth]{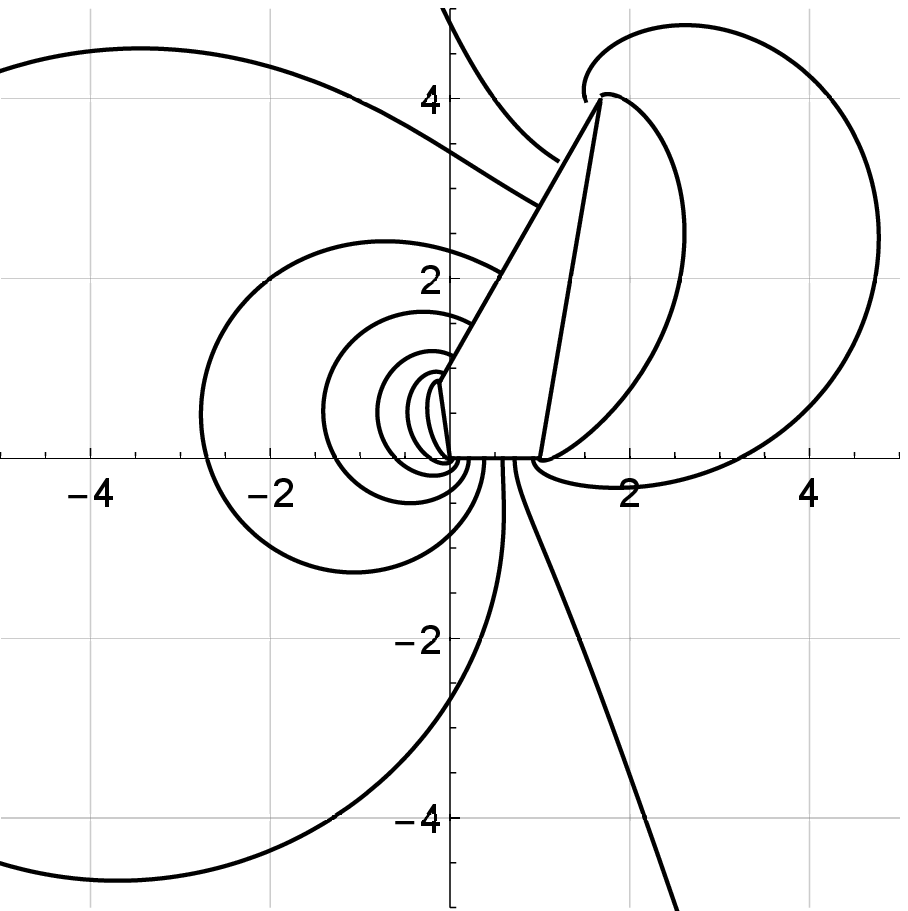}}
\end{minipage}
\hfill
\begin{minipage}[h]{0.49\linewidth}\vspace{3mm}
\center{\includegraphics[width=0.9\linewidth]{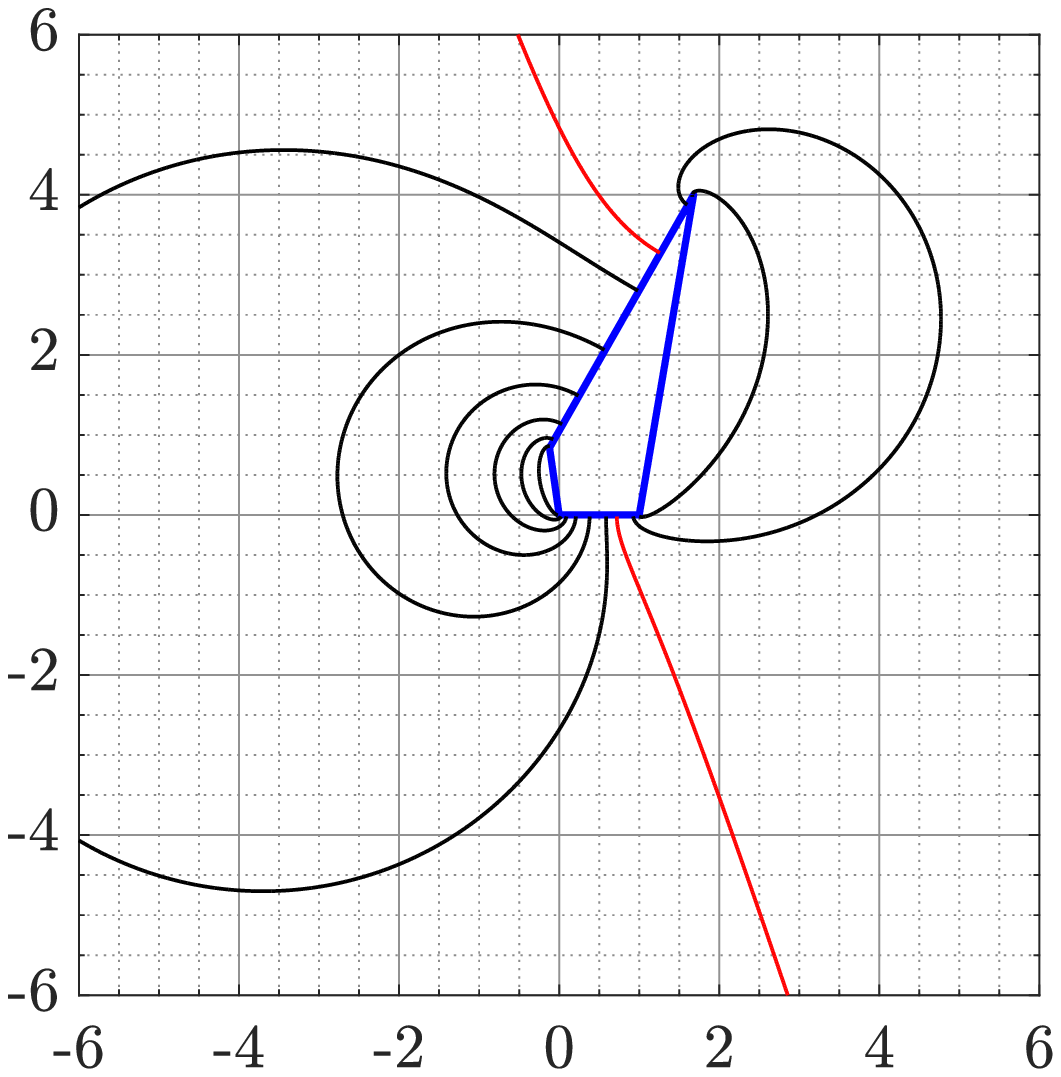}}
\end{minipage}
\caption{The exterior polygon \ref{ep2} and the level curves of its potential function. The figure on the left/right side is produced with the Mathematica/MATLAB codes.}
\label{extp2}
\end{figure}

\begin{nonsec}{\bf Another example.} \label{ep2}
 The polygon shown in Fig.\ref{extp2} has vertices
$1, 0, -3/25 + \i\, 21/25, 42/25 + \i\,4$. The levels are:
$ 0.1, 0.2, u(\infty)=0.334052..., 0.4, 0.5, 0.6, 0.7, 0.8$ and $0.9.$
\end{nonsec}


\section{Exterior of the unit disk}

We consider here the case of a circular quadrilateral determined by the unit disk and four points
on the unit circle. Let $Q$ be the quadrilateral  which is the exterior of the unit disk with
vertices $e^{\i\beta}$, $e^{\i\alpha}$, $e^{-\i\alpha}$, $e^{-\i\beta}$, $0<\alpha<\beta<\pi$.
The M\"obius  transformation
$$
w=T(z)=\frac{z-a}{1-az}, \quad a=\frac{1-\sqrt{\tan \frac{\alpha}{2}\tan \frac{\beta}{2}}}{1+\sqrt{\tan \frac{\alpha}{2}\tan \frac{\beta}{2}}}\,,
$$
maps it onto the quadrilateral $Q_1$ which is the exterior of the unit disk with vertices
$-e^{-\i\gamma}$, $e^{\i\gamma}$, $e^{-\i\gamma}$, $-e^{-\i\gamma}$, where
$$
\cos\gamma=\frac{\sin\frac{\beta-\alpha}{2}}{\sin\frac{\beta+\alpha}{2}}\,,\quad 0<\gamma<\pi/2.
$$
The M\"obius transformation
$$
\omega=S(w)=\frac{1}{\sqrt{\lambda}}\,\frac{w-\i}{1-\i w}
$$
maps $Q_1$  onto the quadrilateral $Q_2$ which is the upper half plane with vertices
$-1/\lambda$, $-1$, $1$, $1/\lambda$, where
$$
\lambda=1+2\tan^2\gamma-2\tan\gamma\sqrt{1+\tan^2\gamma}.
$$
At last,
$$
\zeta=g(\omega)=\frac{F(\lambda,\arcsin \omega)}{2\K(\lambda)}\,+0.5
$$
maps  $Q_2$ onto $Q_3$ which is the rectangle $[0,1]\times [0,h]$ with vertices $0$, $1$, $1+\i h$, $\i h$,
where $h$ is the conformal modulus of $Q$.

The function
$$
u(z)=\Re g\circ S\circ T(z)
$$
is the potential function for $Q$. We note that
$g\circ S\circ T(\infty)=g\circ S(-1/a)=g\left(-\frac{1}{\sqrt{\lambda}}\,\frac{1+\i a}{a+\i}\right)$,
therefore, for $u_0=u(\infty)$
we have
$$
u_0=\Re g\left(-\frac{1}{\sqrt{\lambda}}\,\frac{1+\i a}{a+\i}\right).
$$
If we put $u_0=u(\infty)$, then the parametric representation of the level curve $u=u_0$ takes the form
$$
z=T^{-1}\circ S^{-1}\circ\, \sn((2u_0-1)\K(\lambda)+\i\xi,\lambda),\quad 0\le\xi  \le \K'(\lambda).
$$
or
$$
z=\frac{\sqrt{\lambda}(a-\i)\,\sn((2u_0-1)\K(\lambda)+\i\xi,\lambda)+(1-\i a)}{\sqrt{\lambda}(1-\i a)\,\sn((2u_0-1)\K(\lambda)+ i\xi,\lambda)+(a- \i)}\,,\quad 0\le\xi  \le  \K'(\lambda).
$$
Here $\sn(\,\cdot\,,\lambda)$\, is the Jacobi sine.

The level curve $\{u=u_0\}$ is symmetric with respect to the real axis. If we consider its upper half,
then, as numerical calculations show, we can see that it  is very close to a (rectilinear) ray and
is essentially distinct from it only in a small neighborhood of its endpoint lying on the unit circle.

\begin{center}
\begin{figure}[H] 
\includegraphics[width=7cm]{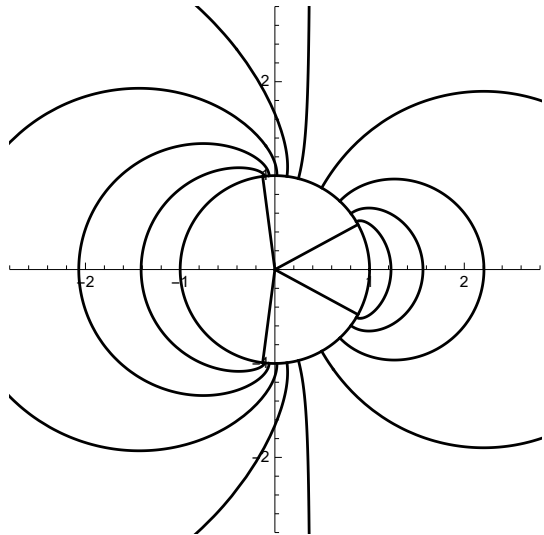}\\[5mm]
\caption{The level curves of the potential function for the mixed
Dirichlet-Neumann problem corresponding to the levels
$0.2j\cdot u_0$, $1\le j\le 5$; $u_0 + 0.2k(1 - u_0)$, $1\le k\le 4$;
$\alpha=0.5$, $\beta=1.7$. The quadrilateral is the exterior of the unit
disk and its vertices are indicated by the endpoints
of the radii.}\label{infDN3}
\end{figure}
\end{center}

The Fig.\ref{infDN3} was produced by the following Mathematica script.
\small
\begin{verbatim}
alpha = 0.5; beta = 1.7; a = (1 - Sqrt[Tan[alpha/2]Tan[beta/2]])/
(1 + Sqrt[Tan[alpha/2] Tan[beta/2]]);
gamma = ArcCos[Sin[(beta - alpha)/2]/Sin[(beta + alpha)/2]];
lambda = 1 + 2 (Tan[gamma])^2 - 2*Tan[gamma]* Sqrt[1 + (Tan[gamma])^2];
zeta0 = (-1/Sqrt[lambda]) (1 + I*a)/(a + I);
u0 = Re[EllipticF[ArcSin[zeta0], lambda^2]/(2*EllipticK[lambda^2])] + 0.5;
u1 = 0.2*u0; u2 = 0.4*u0; u3 = 0.6*u0; u4 = 0.8*u0; u5 = u0 + 0.2 (1 - u0);
u6 = u0 + 0.4 (1 - u0); u7 = u0 + 0.6 (1 - u0); u8 = u0 + 0.8 (1 - u0);
t0 = EllipticK[1 - lambda^2]/EllipticK[lambda^2];
H[xi_, u_] = (Sqrt[lambda] (a - I) JacobiSN[(2*u - 1 + I*xi) (EllipticK[lambda^2]),
lambda^2] + (1 - a*I))/(Sqrt[lambda] (1 - a*I) JacobiSN[(2*u - 1 + I*xi)
(EllipticK[lambda^2]), lambda^2] + (a - I));
figInf3 = ParametricPlot[{{Re[H[t, u0]],Im[H[t, u0]]}, {Re[H[t, u1]], Im[H[t, u1]]},
{Re[H[t, u2]],Im[H[t, u2]]}, {Re[H[t, u3]], Im[H[t, u3]]}, {Re[H[t, u4]],
Im[H[t, u4]]}, {Re[H[t, u5]], Im[H[t, u5]]}, {Re[H[t, u6]],Im[H[t, u6]]},
{Re[H[t, u7]], Im[H[t, u7]]}, {Re[H[t, u8]],Im[H[t, u8]]}, {Cos[2*Pi*t/t0],
Sin[2*Pi*t/t0]}, {(t/t0) Cos[alpha], (t/t0) Sin[alpha]}, {(t/t0) Cos[beta],
(t/t0) Sin[beta]}, {(t/t0) Cos[alpha], -(t/t0) Sin[alpha]}, {(t/t0) Cos[beta],
-(t/t0) Sin[beta]}}, {t, 0, t0}, PlotRange -> 2.8, PlotStyle -> Black,
AxesStyle -> Black]
\end{verbatim}

\normalsize
The radii in the figure show the location of the vertices of the quadrilateral.

\section{Some numerical results}

Now we will give  values of $u(\infty)$ for some unbounded polygonal
quadrilaterals $$Q=(D;z_1,z_2,z_3,z_4)$$ with vertices  $z_1=1$, $z_2=0$, $z_3=B$,
$z_4=A$ and angles $\pi(1+\alpha_k)$, $1\le k\le 4$. To calculate $u(\infty)$,
we use the Mathematica module ExtMod[$B,A,n,wp$]
from \cite[Appendix~A]{nsv}
written with the help of the Wolfram Mathematica system.
This function returns the values of
$\text{Mod}(Q)$, $\alpha_1$, $\alpha_2$, $\alpha_3$, $\alpha_4$, $t$,
and $z_0$.
For example, if $A= 7 + 5\, \i$, $B= -1 + 2\, \i$ (see Table~\ref{tab3}, Line~2),
then we use the module  ExtMod[$B,A,n,wp$] as follows:
\bigskip

\small
\begin{verbatim}
uInfty[A_, B_] := Module[{v = ExtMod[B,A,2,16], r, t, z0}, t = v[[6]]; z0 = v[[7]];
r = 1/Sqrt[t]; Re[EllipticF[ArcSin[Sqrt[z0]], r^2]/EllipticK[r^2]]];
DecimalForm[uInfty[7 + 5 I,-1 + 2 I ],16]
\end{verbatim}

\normalsize

\medskip

\noindent
In the above codelines we have used the notation of  \eqref{uinf}, to make
the application of \newline ExtMod[$B,A,n,wp$] transparent for the readers.
Here the parameter $wp$ stands for the accuracy goal, and the above
call {\tt  ExtMod[B,A,2,16]} means that the desired accuracy is $10^{-16}.$
\bigskip

\begin{table}[ht]
\caption{The values of $u(\infty)$ for some polygonal quadrilaterals
computed using the formula~\eqref{uinf} and using the PlgCirMap
MATLAB toolbox~\cite{pcm}. The computed values  agree with $9$
decimal places or more.
}
\centering
\begin{tabular}{|l|l|l|l|}
 \hline
\hspace{3mm} $A$ & \hspace{3mm} $B$ & Using the Formula~\eqref{uinf} & PlgCirMap toolbox method \\
 \hline
 7 + 5 I &$ -1 + 2$ I  &0.3782951219491777 &0.378295121963035\\   \hline
 8 + 3 I &  $-1$ + I   &0.3507184214435048 &0.350718421549051\\  \hline
 5 + 5 I &  $-3$ + I   &0.4209495357540314 &0.420949535761708\\   \hline
  7 + 4 I & $-3$ + 3 I &0.4473431220217027 &0.447343122027837\\    \hline
   5 + 5 I& $-1$ + 2 I &0.3916188047098933 &0.391618804701488\\   \hline
  7 + 5 I&\hspace{8mm}I&0.3172197705784933 &0.317219770131948\\   \hline
  7 + 3 I&\hspace{2mm}1 + 2 I &0.3917841755037506 &0.391784175504020\\  \hline
   4 + 5 I &  $-2$ + I    &0.3960930352825737 &0.396093035293408\\   \hline
     1 + I& \hspace{8mm}I &0.5000000000000000 &0.500000000000000\\       \hline
\end{tabular}\label{tab3}
\end{table}

\end{document}